\documentclass[a4paper,12pt, draft]{article}
\usepackage{a4wide}
\usepackage{amsmath,amsthm,amssymb}
\usepackage[mathscr]{eucal}
\usepackage{enumerate}
\theoremstyle{plain}
\newtheorem{proposition}{Proposition}
\newtheorem{theorem}{Theorem}
\newtheorem{lemma}{Lemma}
\theoremstyle{definition}

\newtheorem{corollary}{Corollary}
\newcommand{\IS}{\mathcal {IS}_n}
\newcommand{\ISA}{(\mathcal {IS}_n, *_a)}
\newcommand{\ISB}{(\mathcal {IS}_n, *_b)}
\newcommand{\ISZ}{(\mathcal {IS}_n, *)}
\renewcommand{\phi}{\varphi}
\renewcommand{\epsilon}{\varepsilon}
\newcommand{\im}{\operatorname{im}}
\newcommand{\ran}{\operatorname{im}}

\newcommand{\rank}{\operatorname{rank}}
\newcommand{\dom}{\operatorname{dom}}
\newcommand{\Aut}{\operatorname{Aut}}
\newcommand{\Ann}{\operatorname{Ann}}
\begin{document}
\author{G.M. Kudryavtseva, G.Y. Tsyaputa}
\title{The automorphism group of the sandwich inverse symmetric semigroup}
\date{}
\maketitle
\begin{abstract}
The structure of the automorphism group of the sandwich semigroup
$\IS$ is described in terms of standard group constructions.
\end{abstract}

\section{Introduction}
Let $X$ and $Y$ be two nonempty sets, $S$ be the set of all maps
from $X$ to $Y$. Fix some $\alpha : Y\rightarrow X$ and define the
multiplication in $S$ in the following way: $\phi\circ \psi =
\phi\alpha\psi$ (the composition of the maps is from the left to
the right). The action defined by this rule is associative. Ljapin
in ~\cite[p.353] has set the problem of investigation of the
properties of this semigroup depending on the restrictions on $S$
and $\alpha$.

 Magill in \cite{Mag} has considered this problem in the case when $X$ and $Y$ are
 topological spaces and the maps are continuous.
 In particular, under the assumption that $\alpha$ is onto he has described the automorphisms
 of such a semigroup and has determined the isomorphism criterion of two such semigroups.

Sullivan in \cite{Sul} has proved that if $|Y|\leq|X|$ then
Ljapin's semigroup embeds into transformation semigroup on the set
$X\cup \{a\}$, $a\notin X$.

The important case is if $X=Y$. In this case one considers $T_X$,
the transformation semigroup on $X$, $\alpha\in T_X$. Symons in
\cite{Sym} has established the isomorphism criterion for such
semigroups and has described their automorphism groups.

The problem of the description of automorphism group is one of the
most important in the process of study of a certain algebraic
structure. In particular, much attention has been devoted to the
description of automorphism groups of different semigroups (see,
for example, \cite{5}, \cite{6} and the references therein).

In the present paper we study the structure of the automorphism
group of the sandwich finite inverse symmetric semigroup $\IS$ of
all partial injective transformations of set $N=\{1,\dots,n\}$
with the sandwich element $a\in\IS$ (denote it by $\ISA$). The
main result is analogous to the one on the automorphism group of
the sandwich semigroup of all transformations obtained by Symons
in \cite{Sym}.

\section{The main theorem}

For an element $a\in\IS$ denote by $\dom(a)$ the domain of $a$,
and by $\im(a)$ the image of $a$. The cardinality of the set
$\im(a)$ is called the \textit{rank} of $a$ and is denoted by
$\rank(a)$. If $\dom(a)=\{x_1,\dots,x_k\}$ and $a(x_i)=y_i$ for
all $1\leq i\leq k$ then we write $a=\left(
 \begin{array}{ccccccccc} x_1 & \dots & x_k \\ y_1 & \dots & y_k
 \end{array}\right)$.

For a subset $A\subseteq N$ denote by $\overline{A}$ the
complement ($N\setminus A$) of $A$. If $a\in \IS$ and $M\subseteq
N$ then $a|_M$ denotes the restriction of $a$ to $M$.

Denote by $E(S)$ the set of idempotents of the semigroup $S$.
Recall (see \cite{id}) that the relation $f\leq h \Leftrightarrow
fh=hf=f$ defines the \textit{natural partial order} on $E(S)$.

Finally for arbitrary $K$ denote by $\mathcal{S}(K)$ the full
symmetric group on $K$. Other notions, which are used in the paper
without definitions, can be found, for example, in \cite{id}.

We will use the following fact first proved in \cite{8}.
\begin{theorem}
Semigroups $\ISA$ and $\ISB$ are isomorphic if and only if
$\rank(a)=\rank(b)$.
\end{theorem}
Hence, up to isomorphism one can study only sandwich semigroups
$\ISZ$ with sandwich elements $e$ such that $e\in E(\IS)$. From
now on we fix idempotent $e$ such that $e$ acts identically on set
$A=\{1,2,\dots, k\}$, $A\subseteq N$, and denote $(\IS, *_e)$ just
by $\ISZ$.

\begin{proposition}\label{lema2.1}
The lattice of idempotents $E\ISZ$ is isomorphic to the lattice
$\mathcal{B}(A)$ of the subsets of $A$.
\end{proposition}
\begin{proof}
In \cite{8} it is proved that element $a\in\IS$ is idempotent in
$\ISZ$ if and only if $a$ is the idempotent in $\IS$ and
$\dom(a)\subseteq A$. Hence the map $f\mapsto\dom(f)$ defines the
isomorphism between the partially ordered set $E\ISZ$ and the
lattice $\mathcal{B}(A)$.
\end{proof}
\begin{corollary}
The sandwich idempotent $e$ is the maximal element of the lattice
of idempotents in $\ISZ$. In particular, for arbitrary
automorphism $\phi\in \Aut\ISZ$ we have $\phi(e)=e$.
\end{corollary}
Let $S$ be a semigroup, $ a\in S$. Element $a$ is called
\textit{decomposable} in $S$ provided that $a=bc$ for some $b,c\in
S$, otherwise $a$ is called \textit{indecomposable}.

Let $\mathrm{P}$ be the set of all decomposable elements in
$\ISZ$, and $\mathrm{Q}$ be the set of all indecomposable
elements. For every $a\in \mathrm{Q}$ consider the sets
\begin{align*}
M_1 (a)&=\{x\in\dom(a) \;|\; x\in A, \;a(x)\in A\},\\
M_2 (a)&= \{x\in\dom(a) \;| \;x\in A, \;a(x)\in \overline{A}\}, \\
M_3 (a)&=\{x\in\dom(a)\; | \;x\in \overline{A},\; a(x)\in A\},\\
M_4 (a)&=\{x\in\dom(a) \;|\; x\in \overline{A}, \;a(x)\in
\overline{A}\}.
\end{align*}
Obviously $\dom(a)=M_1 (a)\cup M_2 (a)\cup M_3 (a)\cup M_4 (a)$.
For $a, b \in \mathrm{Q}$ set $a \sim b$ provided that $M_i
(a)=M_i (b)$, $1\leq i\leq 3$, and $a(x)=b(x)$ for all $x\in M_1
(a)\cup M_2 (a)\cup M_3 (a)$. If either $a$ or $b$ does not belong
to $\mathrm{Q}$ then set $a\sim b$ provided that $a=b$.
\begin{lemma}\label{lema2.6}
The element $a\in\IS$ is decomposable in $\ISZ$ if and only if
$\rank(a)\leq k$.
\end{lemma}
\begin{proof}
Let $a$ be decomposable in $\ISZ$. Then $a=b * c$ for some
$b,c\in\IS$, therefore $\rank(a)=\rank(bec)\leq\rank(e)=k$.

Conversely, let $\rank(a)=l\leq k$ and $a=\left(
\begin{array}{lccccr} x_{1} &\dots& x_{l}\\
 y_{1} & \dots & y_{l} \end{array}
 \right)$. The direct calculation shows that $a=b *c$, where
 $b=\left(
\begin{array}{lccccr} x_{1} &\dots& x_{l}\\
 1 & \dots & l \end{array}
 \right)$, $c=\left(
\begin{array}{lccccr} 1 &\dots& l\\
 y_{1} & \dots & y_{l} \end{array}
 \right)$.
\end{proof}
\begin{lemma}\label{lema2.7}
The relation $\sim$ is a congruence on $\ISZ$.
\end{lemma}
\begin{proof}
Let $a\sim a_1$ and $b\sim b_1$. Show that $a *b \sim a_1 *b_1$.
Since $\rank(a *b)\leq k$ and $\rank(a_1 *b_1)\leq k$, using Lemma
\ref{lema2.6} we get $a *b,\;a_1 *b_1 \in \mathrm{P}$. Hence by
the definition of $\sim$ we have to prove that $a *b=a_1 *b_1$. In
fact, $a(x)\in A$, $a_1 (x)=a(x)$, and $b(a(x))=b_1 (a(x))$ for
arbitrary $x\in\dom(a *b)$. Hence $\left(a
*b\right)(x)=b(a(x))=b_1 (a(x))=b_1 (a_1 (x))=\left(a_1
*b_1\right)(x)$.

Analogously, for every $x\in\dom(a_1 *b_1)$ we obtain $\left(a_1
*b_1\right)(x)=\left(a *b\right)(x)$.
\end{proof}
\begin{theorem}\label{theor2.7}
Let $\ISZ$ be the sandwich semigroup $\IS$ with sandwich
idempotent $e$ and let $\dom(e)=A$. Let further $P_1,\dots,P_m$ be
the list of equivalence classes of the relation $\sim$ on $\IS$.
Then
\begin{displaymath}
\Aut\ISZ=\left(\overset{m}{\underset{i=1}{\oplus}}\mathcal{S}(P_i)\right)\leftthreetimes\Bigl(
\mathcal{S}(A)\times\mathcal{S}(\overline{A})\times\mathcal{S}(\overline{A})\Bigr).
\end{displaymath}
\end{theorem}
We split the proof into several lemmas.

\begin{lemma}\label{lema2.8}
$\left(\overset{m}{\underset{i=1}{\oplus}}\mathcal{S}(P_i)\right)<\Aut\ISZ$.
\end{lemma}
\begin{proof} The statement follows from the definition
of $\sim$ and Lemma \ref{lema2.7}.
\end{proof}

Denote by $\mathcal{G}$ the group $\left(
\mathcal{S}(A)\times\mathcal{S}(\overline{A})\times\mathcal{S}(\overline{A})\right)$.
For arbitrary $g\in\mathcal{S}(A)$ and
 $h\in\mathcal{S}(\overline{A})$ we denote by
 $_{g}h$ the element $(g, h)$ in $\left(
\mathcal{S}(A)\times\mathcal{S}(\overline{A})\right)$ . It is
obvious that $_{g}h$ can be regarded as the element in $\IS$ via
the identification with its image under the natural embedding:
$\left(\mathcal{S}(A)\times\mathcal{S}(\overline{A})\right)\hookrightarrow\IS$.
.

\begin{lemma}\label{lema2.9}
For arbitrary $(g, h_1, h_2)\in \mathcal{G}$ the map
\begin{displaymath}
\tau_{(g, h_1, h_2)} \; :\; a\;\mapsto\; _{g}h_{1}^{-1} a
\,_{g}h_{2}
\end{displaymath}
is an automorphism of $\ISZ$, moreover the map
 $(g, h_1 ,h_2)\overset{\phi}{\longrightarrow} \tau_{(g, h_1, h_2)}$
is the monomorphism from $\mathcal{G}$ to $\Aut\ISZ$.
\end{lemma}

\begin{proof}
Set $\tau=\tau_{(g, h_1, h_2)}$. It is enough to show that $\tau(a
*b)=\tau(a)*\tau(b)$ for arbitrary $a, b\in \IS$. By the
construction of $ _{g}h_{1}$ and $_{g}h_{2}$ we get
$\;_{g}h_{1}^{-1} * \,_{g}h_{2}= \; _{g}h_{1}^{-1} e
\,_{g}h_{2}=e$, hence
\begin{displaymath}
\tau(a *b)=\; _{g}h_{1}^{-1} a *b \,_{g}h_{2}=\; _{g}h_{1}^{-1} a
e b \,_{g}h_{2}=\; _{g}h_{1}^{-1} a\,_{g}h_{2}*\; _{g}h_{1}^{-1} b
\,_{g}h_{2}=\tau(a)*\tau(b).
\end{displaymath}
\end{proof}

We identify $\mathcal{G}$ with  $\phi(\mathcal{G})$ and consider
$\mathcal{G}$ as a subgroup of $\Aut\ISZ$.

\begin{lemma}\label{lema2.10}
\begin{displaymath}
\Bigl(
\mathcal{S}(A)\times\mathcal{S}(\overline{A})\times\mathcal{S}(\overline{A})\Bigr)
\bigcap
\left(\overset{m}{\underset{i=1}{\oplus}}\mathcal{S}(P_i)\right)=\epsilon,
\end{displaymath}
 where $\epsilon$ is the identical automorphism of $\ISZ$.
\end{lemma}

\begin{proof}
Let $\pi\in \Bigl(
\mathcal{S}(A)\times\mathcal{S}(\overline{A})\times\mathcal{S}(\overline{A})\Bigr)
\bigcap
\left(\overset{m}{\underset{i=1}{\oplus}}\mathcal{S}(P_i)\right)$.
Take arbitrary $t,s\in N$ and $a=\left(
\begin{array}{lccr} t
\\ s \end{array} \right)\in \IS$. Then $a\in\mathrm{P}$. This and
$\pi\in \overset{m}{\underset{i=1}{\oplus}}\mathcal{S}(P_i)$ imply
$\pi(a)=a$. Since $\pi\in
\mathcal{S}(A)\times\mathcal{S}(\overline{A})\times\mathcal{S}(\overline{A})$
it follows that $\pi=\tau_{(g,h_1, h_2)}$ for some
$g\in\mathcal{S}(A)$, $h_1 ,h_2 \in\mathcal{S}(\overline{A})$.
Then $\pi(a)=\pi\left(
\begin{array}{lccr} t
\\ s \end{array} \right) =\left( \begin{array}{lccr} \; _{g}h_{1}
(t)
\\ \; _{g}h_{2} (s) \end{array} \right)$. Hence, $\;_{g}h_{1}(t)=t$
 and $\; _{g}h_{2}(s)=s$, therefore $\; _{g}h_{1}=\; _{g}h_{2}=\epsilon$.

 Hence $\pi(a)=\; _{g}h_{1}^{-1}a\;
 _{g}h_{2}=a$ for arbitrary $a\in\IS$, and $\pi=\epsilon$.
\end{proof}

\begin{lemma}\label{lema2.11}
Let $\pi\in\Aut\ISZ$ and $a\in\mathrm{P}$. Then
$\rank(\pi(a))=\rank(a)$.
\end{lemma}
\begin{proof}
Proposition~\ref{lema2.1} implies that $\rank(\pi(f))=\rank(f)$
for any idempotent $f\in E\ISZ$. Let $a\in\mathrm{P}$. For
arbitrary idempotent $f\in E\ISZ$ we have
\begin{displaymath} \IS *f *\IS=\{a\in\IS
\,:\,\rank(a)\leq\rank(f)\}.
\end{displaymath}
Hence $a$ is contained in the main two-sided ideal generated by
$f$ if and only if $\rank(f)\geq\rank(a)$. This yields
\begin{displaymath}
\rank(a)=\min\{k\,:\,a\in \IS *f *\IS \text{ and } f \text{ is
idempotent of rank } k\}.
\end{displaymath}
It follows that
\begin{displaymath}
\rank(\pi(a))=\min\{k\,:\,\pi(a)\in \IS *q *\IS \text{ and } q
\text{ is idempotent of rank  } k\}=\rank(a).
\end{displaymath}
\end{proof}
Consider the sets:
\begin{displaymath} \Ann_L \ISZ=\{x\in\IS : x*a=0 \text{ for all } a\in
\IS\},
\end{displaymath}
\begin{displaymath} \Ann_R \ISZ=\{x\in\IS : a*x=0 \text{ for all }
a\in \IS\},
\end{displaymath}
\begin{displaymath}
\text{and } \Ann\ISZ=\Ann_L \ISZ\cap\Ann_R \ISZ.
\end{displaymath} The following lemma is obvious.
\begin{lemma}\label{ann}
\begin{enumerate}[1)]
\item $x\in\Ann_L \ISZ$ if and only if $\ran(x)\in\overline{A}$,
\item $x\in\Ann_R \ISZ$ if and only if $\dom(x)\in\overline{A}$,
\item $x\in\Ann\ISZ$ if and only if $\ran(x)\in\overline{A}$ and
$\dom(x)\in\overline{A}$.
\end{enumerate}
\end{lemma}
\begin{lemma}\label{lema2.12}
$\Aut\ISZ=\Bigl(
\mathcal{S}(A)\times\mathcal{S}(\overline{A})\times\mathcal{S}(\overline{A})\Bigr)
\cdot
\left(\overset{m}{\underset{i=1}{\oplus}}\mathcal{S}(P_i)\right)$.
\end{lemma}
\begin{proof} We show that for arbitrary $\sigma\in\Aut\ISZ$ there
exist $\tau=(g, h_1, h_2)\in\left(
\mathcal{S}(A)\times\mathcal{S}(\overline{A})\times\mathcal{S}(\overline{A})\right)$
and
$\pi\in\left(\overset{m}{\underset{i=1}{\oplus}}\mathcal{S}(P_i)\right)$,
such that $\sigma=\tau\pi$. For every $i\in A$ let
$f_{\{i\}}=\left(
\begin{array}{lccr} i
\\ i \end{array} \right)$. Proposition \ref{lema2.1}
implies that $\sigma$ induces some permutation $\widetilde{g}$ on
the set $\{f_{\{i\}}, i\in A\}$, which in its turn, induces the
permutation of the indices $g\in\mathcal{S}(A)$.

Lemmas \ref{lema2.11} and \ref{ann} imply that $\left(
\begin{array}{lccr} j \\ j
\end{array} \right)\in\Ann\ISZ$ for arbitrary $j\in
\overline{A}$, moreover, using the equality $\sigma(\Ann\ISZ)=
\Ann\ISZ$, we get $\sigma\left(\begin{array}{lccr} j
\\ j \end{array} \right)=\left( \begin{array}{lccr} s
\\ t \end{array} \right)$, where $s,t\in \overline{A}$.

Set $h_1 (j)=s$, $h_2 (j)=t$. Verify that the maps $h_1 ,\, h_2
\,:\, \overline{A}\rightarrow \overline{A}$ are bijections. Assume
that $h_1 (s)=h_1 (l)$, $k+1\leq s<l\leq n$. This means that
\begin{displaymath}\sigma\left(
\begin{array}{lccr} s
\\ s \end{array} \right)=\left( \begin{array}{lccr} h_1 (s)
\\ h_2 (s) \end{array} \right) \text{ and } \sigma\left( \begin{array}{lccr}
l
\\ l \end{array} \right)=\left( \begin{array}{lccr} h_1 (s)
\\ h_2 (l) \end{array} \right).
\end{displaymath}
However in this case
\begin{displaymath}\left( \begin{array}{lccr} s
\\ s \end{array} \right)=\sigma^{-1}\left( \begin{array}{lccr} h_1 (s)
\\ 1 \end{array} \right)*\sigma^{-1}\left( \begin{array}{lccr} 1
\\ h_2 (s) \end{array} \right),
\end{displaymath}
\begin{displaymath}
 \left( \begin{array}{lccr} l
\\ l \end{array} \right)=\sigma^{-1}\left( \begin{array}{lccr} h_1 (s)
\\ 1 \end{array} \right)* \sigma^{-1}\left( \begin{array}{lccr} 1
\\ h_2 (l) \end{array} \right),
\end{displaymath}
 hence, $s=l$ which contradicts the choice of
$s$ and $l$.

Now we prove that $\sigma\left(
\begin{array}{lccr} s
\\ t \end{array} \right)=\left( \begin{array}{lccr} \;_{g}h_1 (s)
\\ \;_{g}h_2(t) \end{array} \right)$, $s,t\in N$. Let $\sigma\left( \begin{array}{lccr} s
\\ t \end{array} \right)=\left( \begin{array}{lccr} u
\\ v \end{array} \right)$. If $s,t\in A$ then $f_{s}*\left( \begin{array}{lccr} s
\\ t \end{array} \right)=\left( \begin{array}{lccr} s
\\ t \end{array} \right)$. Taking the $\sigma$ pre-images
we get : $f_{g(s)}*\left( \begin{array}{lccr} u
\\ v \end{array} \right)=\left( \begin{array}{lccr} u
\\ v \end{array} \right)$, hence $u=g (s)=\;_{g}h_1 (s)$. In the same vein we get $\left( \begin{array}{lccr} u
\\ v \end{array} \right)*\epsilon_{\{g(t)\}}=\left( \begin{array}{lccr}
u
\\ v \end{array} \right)$ from the equality $\left( \begin{array}{lccr} s
\\ t \end{array} \right)*\epsilon_{\{t\}}=\left( \begin{array}{lccr} s
\\t \end{array} \right)$, and $v=g(t)=\;_{g}h_2 (t)$.

Now let $s\in \overline{A}$, $t\in A$. Since $\left(
\begin{array}{lccr} s
\\ t \end{array} \right) \in\Ann_R \ISZ$ then $\sigma\left( \begin{array}{lccr} s
\\ t \end{array} \right)=\left( \begin{array}{lccr} u
\\ g(t) \end{array} \right)\in \Ann_R \ISZ$, and $u\in \overline{A}$.
Analogously since $\left( \begin{array}{lccr} t
\\ s \end{array} \right)\in\Ann_L \ISZ$, then $\sigma\left( \begin{array}{lccr}
t
\\ s \end{array} \right)=\left( \begin{array}{lccr} g(t)
\\ v \end{array} \right)\in \Ann_L \ISZ$, and $v\in \overline{A}$.
Applying $\sigma$ to the both sides of the equality $\left(
\begin{array}{lccr} s
\\ s \end{array} \right)=\left( \begin{array}{lccr} s
\\ t \end{array} \right)*\left( \begin{array}{lccr} t
\\ s \end{array} \right)$ we get $\left( \begin{array}{lccr} h_1(s)
\\ h_2 (s) \end{array} \right)=\left( \begin{array}{lccr} u
\\ g(t) \end{array} \right)*\left( \begin{array}{lccr} g(t)
\\ v \end{array} \right)$. Hence $u=h_1 (s)$, $v=h_2 (s)$.

The case when $s\in A$, $t\in \overline{A}$ is treated in the same
way.

Finally, let $s\in \overline{A}$, $t\in \overline{A}$. The
equality
\begin{align*}\left( \begin{array}{lccr} u
\\ v \end{array} \right)=\tau\left( \begin{array}{lccr} s
\\ t \end{array} \right)=\tau\left(\left( \begin{array}{lccr} s
\\ 1 \end{array} \right)*\left( \begin{array}{lccr} 1
\\ t \end{array} \right)\right)=\\=\tau\left( \begin{array}{lccr} s
\\ 1 \end{array} \right)*\tau\left( \begin{array}{lccr} 1
\\ t \end{array} \right)=\left( \begin{array}{lccr} h_1(s)
\\ g(1) \end{array} \right)*\left( \begin{array}{lccr}
g(1)
\\ h_2 (t) \end{array} \right)
\end{align*}
implies that $u=h_1 (s)$, $v=h_2 (t)$.

Let $a\in\IS$. We show that
\begin{equation}\label{eq2.11}
 M_1 (\sigma(a))=g(M_1 (a)) \text{ and } (\sigma(a))(g(x))=g(a(x)),
 \text{ for all } x\in M_1 (a),
\end{equation}
\begin{equation}\label{eq2.12}
 M_2 (\sigma(a))=g(M_2 (a)) \text{ and } (\sigma(a))(g(x))=h_2 (a(x)),
 \text{ for all } x\in M_2 (a),
\end{equation}
\begin{equation}\label{eq2.13}
 M_3 (\sigma(a))=h_1 (M_3 (a)) \text{ and } (\sigma(a))(h_1 (x))=g(a(x)),
 \text{ for all } x\in M_3 (a).
\end{equation}

Prove, for example, (\ref{eq2.12}), then statements (\ref{eq2.11})
and (\ref{eq2.13}) can be proved analogously.

Let $x\in M_2(a)$. Then $x\in A$, $a(x)\in \overline{A}$. From the
equality $\left( \begin{array}{lccr} x
\\ x \end{array} \right)*a=\left( \begin{array}{lccr} x
\\ a(x) \end{array} \right)$,
using the arguments from the previous paragraph, we get
\begin{displaymath}\left(
\begin{array}{lccr} g(x) \\ g(x) \end{array} \right)*\sigma(a)=
\sigma\left(\left( \begin{array}{lccr} x
\\ x \end{array} \right)*a\right)=\sigma\left( \begin{array}{lccr}
x
\\ a(x) \end{array} \right)=\left( \begin{array}{lccr} g(x)
\\ h_2(a(x)) \end{array} \right).
\end{displaymath}
Therefore $g(x)\in\dom(\sigma(a))$ and $(\sigma(a))(g(x))=h_2
(a(x))$. In particular, $g(x)\in M_2 (\sigma(a))$. It follows that
$g(M_2 (a))\subseteq M_2(\sigma(a))$ and $(\sigma(a))(g(x))=h_2
(a(x))$ for all $x\in M_2 (a)$.

Conversely, let $x\in M_2 (\sigma(a))$. Denote $y=(\sigma(a))(x)$.
Then $x\in A$, $y\in \overline{A}$. Going to the $\sigma^{-1}$
images in the equality $\left(
\begin{array}{lccr} x
\\ x \end{array} \right)*\sigma(a)=\left( \begin{array}{lccr} x
\\ y \end{array} \right)$, we get $\left(
\begin{array}{lccr} g^{-1}(x)
\\ g^{-1}(x) \end{array} \right)*a=\left( \begin{array}{lccr}
g^{-1}(x)
\\ h_2^{-1}(y) \end{array} \right)$, hence $g^{-1}(x)\in\dom(a)$ and
$a(g^{-1}(x))=h_2^{-1}(y)$, moreover, $h_2^{-1}(y)\in
\overline{A}$. Then $g^{-1}(x)\in M_2 (a)$, and $x\in g(M_2 (a))$.

Now let $a=\left( \begin{array}{lcccccr} x_1 &\dots & x_l
\\ y_1 & \dots & y_l \end{array} \right)\in \mathrm{P}$. Lemma \ref{lema2.6} implies that $l\leq k$.
Then $a=b *c$, where $b= \left(
\begin{array}{lccccr} x_1 & \dots & x_l
\\ 1 & \dots & l \end{array} \right)$, $c=\left( \begin{array}{lcccccr} 1 & \dots &
l
\\ y_1 & \dots & y_l \end{array} \right)$. Therefore $\sigma(a)=\sigma(b)*\sigma(c)$, and
hence
\begin{displaymath}
\dom(\sigma(b))=\{\;_{g}h_1 (x_1),\dots, \,_{g}h_1 (x_l)\},\quad
\dom(\sigma(c))=\{g(1),\dots,g(l)\}.
\end{displaymath}
Moreover, for all $1\leq i\leq l$ we have
 \begin{displaymath}\sigma(b)(\:_{g}h_1
(x_i))=g(i),\quad \sigma(c)(g(i))=\,_{g}h_2 (y_i).
\end{displaymath}
Hence,
\begin{displaymath}
\dom(\sigma(a))=\{\:_{g}h_1 (x_1),\dots,\:_{g}h_1 (x_l)\}
\quad\text{ and }
\end{displaymath}
\begin{displaymath}
\sigma(a)(\:_{g}h_1 (x_i))=\:_{g}h_2 (y_i),\quad 1\leq i\leq l.
\end{displaymath}

It follows from (\ref{eq2.11}), (\ref{eq2.12}), (\ref{eq2.13})
 that
\begin{displaymath}
M_i(\sigma(a))=M_i(\tau(a)) \text{ and } \sigma(a)(x)=\tau(a)(x)
\end{displaymath}
for all $x\in M_i (\sigma(a))$, $1\leq i\leq 3$. Hence $a\sim b$
if and only if $\sigma(a)\sim \sigma(b)$, that is if and only if
$\tau(a)\sim \tau(b)$. So $\sigma$ induces the permutation $\pi\in
\overset{m}{\underset{i=1}{\oplus}}\mathcal{S}(P_i)$ such that
$\sigma=\tau\cdot \pi$.
\end{proof}
\begin{lemma}\label{lema2.13}$\displaystyle\overset{m}{\underset{i=1}{\oplus}}\mathcal{S}(P_i)
\vartriangleleft\Aut\ISZ$.
\end{lemma}
\begin{proof}
It is enough to prove that $\tau^{-1}\pi\tau\in
\overset{m}{\underset{i=1}{\oplus}}\mathcal{S}(P_i)$ for arbitrary
$\tau\in
\mathcal{S}(A)\times\mathcal{S}(\overline{A})\times\mathcal{S}(\overline{A})$,
$\pi\in \overset{m}{\underset{i=1}{\oplus}}\mathcal{S}(P_i)$. we
prove that $x\sim \tau^{-1}\pi \tau (x)$ for arbitrary $x\in N$.
Indeed, $\tau^{-1}\pi (x)\sim \tau^{-1}(x)$. Hence $\tau^{-1}\pi
\tau(x)\sim \tau^{-1} \tau(x)=x$.
\end{proof}
\textbf{The proof of theorem \ref{theor2.7}} follows immediately
from the definition of the semidirect product, and Lemmas
\ref{lema2.10}, \ref{lema2.12} and \ref{lema2.13}.
\begin{flushright}
$\square$
\end{flushright}

\begin{center}
\bf Acknowledgments
\end{center}

 For the first author the research was partially supported by The Swedish
Institute. The authors are thankful to Dr O. Ganyushkin and Dr V.
Mazorchuk for the useful remarks and the discussion over the
research.

\noindent Department of Mechanics and Mathematics,\\ Kiev Taras
Shevchenko University, \\ 64, Volodymyrska st., 01033, Kiev,
Ukraine
\\ e-mail: G.K.: {\em akudr\symbol{64}univ.kiev.ua} \vspace{0.2cm}

$\;\:\:\quad$ G.T.: {\em gtsyaputa\symbol{64}univ.kiev.ua}


\begin{thebibliography}{99}
\bibitem{Ljap} Ljapin Y.S., Semigroups, Moscow, Fizmatgiz, 1960
(Russian).

\bibitem{Mag} Magill Kenneth D., Semigroup structures for families of
functions. II. Continuous functions. // J. Austral. Math. Soc. 7
(1967), 95-107.

\bibitem{Sul} Sullivan R.P., Generalized partial transformation
semigroups. // J. Austral. Math. Soc. 19 (1975), part 4, 470-473

\bibitem{Sym} Symons J.S.V., On a generalization of the transformation
semigroup. // J. Austral. Math. Soc. 19 (1975), 47-61

\bibitem{5} Ganyushkin O.G., Temnikov S.G., Kudryavtseva G.M.
The automorphisms groups of maximal nilpotent subsemigoups of
semigroup $\mathcal{IS}(M)$. // Mathematical Studii 1, vol.13
(2000), 11- 22. (in Ukrainian)

\bibitem{6} Szechtman F. On the Automorphism Group of the Centralizer of an
Idempotent in the Full Transformation Monoid // Semigroup Forum.
v.70 (2004), .

\bibitem{id} Artamonov V.A, Salij V.N., Skornyakov L.A. and others,
General Algebra, Moscow, Nauka, 1991, vol. 1 (Russian).

\bibitem{8}Tsyaputa G.Y., Transformation semigroups with the deformed
multiplication. // Bulletin of the University of Kiev,
Series:Physics $\&$ Mathematics., 3 (2003), 82-88 (in Ukrainian).

\end{thebibliography}
\end{document}